\documentclass[letterpaper, 10pt, conference]{ieeeconf}
\IEEEoverridecommandlockouts

\usepackage{cite}
\usepackage{amsmath,amssymb,amsfonts}
\usepackage{graphicx}
\usepackage{textcomp}
\usepackage{xcolor}
\usepackage{xspace}
\usepackage{booktabs}
\usepackage{algorithm}
\usepackage{algpseudocode}
\usepackage{xurl}
\usepackage{tabularx}
\usepackage{siunitx}
\allowdisplaybreaks
\definecolor{lightgray}{gray}{0.9}
\definecolor{darkergray}{gray}{0.8}
\definecolor{forestgreen}{RGB}{34,139,34}

\usepackage{hyperref}
\usepackage[capitalize]{cleveref}
\hypersetup{
    colorlinks=true,
    citecolor=green,
    filecolor=black,
    linkcolor=red,
    urlcolor=blue
}


\begin{document}

\title{\LARGE \bf
Fast Relax-and-Round Unit Commitment with Transmission Constraints
}

\author{Evan J. R. Brody$^{1,\dagger}$, Charles Foltz$^{1}$, Eve Tsybina$^{1}$, Slaven Pele\v{s}$^{1}$, Shaked Regev$^{1,\dagger}$% <-this % stops a space
\thanks{Notice: This manuscript has been authored by UT-Battelle, LLC, under contract DE-AC05-00OR22725 with the US Department of Energy (DOE). The US government retains and the publisher, by accepting the article for publication, acknowledges that the US government retains a nonexclusive, paid-up, irrevocable, worldwide license to publish or reproduce the published form of this manuscript, or allow others to do so, for US government purposes. DOE will provide public access to these results of federally sponsored research in accordance with the DOE Public Access Plan (\href{https://www.energy.gov/doe-public-access-plan}{https://www.energy.gov/doe-public-access-plan}). Research sponsored by the Laboratory Directed Research and Development Program of Oak Ridge National Laboratory, managed by UT-Battelle, LLC, for the US Department of Energy.}% <-this % stops a space
\thanks{$^{1}$Oak Ridge National Laboratory, 1 Bethel Valley Road, Oak Ridge, Tennessee, USA.
        {\tt\small regevs@ornl.gov}}%
\thanks{$^{\dagger}$Equal Contribution}%
}

\maketitle

\begin{abstract}
Recent developments in the US knowledge economy have created a significant growth in datacenter loads, with two major consequences for power generation. First, to compensate for growth in load, datacenters are encouraged to bring their own generating units. Second, in search of the remaining pools of dispatchable generation, utilities are increasingly turning to subtransmission and distribution level generating assets. Coupled with increasing loads and the resulting tighter grid conditions, both trends are likely to create a need to commit a large number of localized generating units under grid constraints. We propose an extension of Relax-and-Round Unit Commitment (RRUC) that is capable of committing generating units for larger problems faster than conventional methods, while staying within intertemporal and spatial MVA constraints. We demonstrate the performance of RRUC using synthetic congestible test systems ranging from 100 to 20,000 buses. 
RRUC consistently finds low cost solutions, independent of the problem size, and its run time increases sub-quadratically in the number of buses. RRUC can solve the 100 bus system in less than a second and the 20,000 bus system 7 minutes. In contrast, a leading state of the art solver cannot find a feasible solution to the 100 bus system in 15 minutes.
\end{abstract}


\section{Introduction}
\label{sec:introduction}

As the US is transitioning to knowledge economy, it becomes more energy intensive through the emergence of compute-focused loads. According to the S\&P Global Market Intelligence Datacenters and Energy Report, in 2024 as much as 43.8 GW of datacenter capacity were planned for deployment between 2025 and 2030~\cite{spglobal2024}. This growth in load will inevitably put the grid under stress. In addition to the unprecedented total volume connecting to the grid, an individual datacenter represents a large local load: the same 2024 S\&P Global Report~\cite{spglobal2024} identifies 222 planned or under-construction datacenters that are 100 MW or larger, including 42 datacenters that are 1 GW or larger. For comparison, in the 2027/2028 PJM capacity auction, cleared capacity volumes ranged from 928.9 MW to 51.7 GW, with most areas between 1.5 GW and 8.5 GW~\cite{PJM_auction}. The interconnection of a single 100 MW data center can add up to 10\% of total area load, and the connection of a 1 GW data center can almost double the load in some areas. Such proportions will make existing transmission constraints more tight and reveal transmission bottlenecks which were not apparent before.

%Simultaneously, 
As grid conditions are becoming more stressed, the number of generating units in the system is projected to increase. The growth in datacenter construction has caused utilities and transmission operators to encourage datacenters to bring their own generation~\cite{ferc2025, bloomberg_podcasts_data_2025}. \cite{btm_2025} reports that as many as 25\% of datacenters under construction are bringing their own generation, usually natural gas, less often inverter based resources 
%(\textbf{IBR}) - commenting as this is the only use
with storage. Connecting this behind-the-meter generation to the grid would require making generating units larger than 20 MW known and potentially available for grid operations. As a result, utilities and transmission operators are likely to face a surge in the number of generating units that they may have to commit within grid constraints.

Further, as more load enters the grid, utilities are increasingly turning to the remaining grid assets to provide peak generation or ancillary services. Usually, these are smaller units located in sub-transmission and distribution networks and operated by smaller utilities nested under larger Regional Transmission Organizations (\textbf{RTOs}),
or aggregations of battery and thermal storage devices \cite{DOE_distribution}. For the purposes of unit commitment (\textbf{UC}), industrial diesel and gas backup turbines \cite{NREL_distribution} and mobile turbines \cite{mobile_generators} are of particular interest.

The present best UC practices lack the scalability needed to accommodate a large number of generating resources. They use mixed integer programs (\textbf{MIP}s) \cite{achterberg2013} solved with heuristics such as branch-and-bound, whose poor scaling properties make them inefficient for large problems \cite{lawler1966, morrison2016}. Existing UC literature studies relatively small systems: 20-130 units in \cite{falvo2022}, 5-1870 units in \cite{montero2022}, up to 100 units in \cite{hong2021}, 5-2709 units in \cite{vanackooij2018}. %(Note: the original paper referenced by \cite{vanackooij2018} uses 2704 generators \cite{fu2007}). 
The number of units does not increase for studies that specifically suggest methods with higher computational efficiency: we found the number of generating units to be between 130 \cite{kim2019} and 289 \cite{yu2024}. 

Further, the MIP formulation does not handle the smooth nonlinear functions that define transmission constraints. Including transmission constraints into a conventional UC problem requires piecewise-linear approximations of alternating current (\textbf{AC}) power flow ~\cite{nanou2021,tejada2019}. This introduces additional deviations from the original power flow problem.

We attempt to address the current limitations of UC problems by proposing a continuous formulation of UC, with generator and transmission constraints. This study advances our earlier efforts~\cite{regev2026oneperiodUC,regev2026economicUC} by extending the previous formulation to a full nodal transmission constrained AC optimal power flow (\textbf{ACOPF}) model and testing it on synthetic systems ranging from 100 to 20,000 buses.

\cref{sec:equations} introduces our mathematical formulation, input parameters, decision variables, algorithm and implementation. \cref{sec:design} introduces our test case class and explains our modeling assumptions and data sources. \cref{sec:algorithm} introduces our algorithm and gives our implementation details. \cref{sec:results} shows the results of running our algorithm on the test case class. \cref{sec:summary} summarizes our results and provides future research directions. 

\section{Mathematical Formulation}
\label{sec:equations}

Our mathematical formulation of UC extends~\cite{regev2026oneperiodUC, regev2026economicUC} and closely overlaps with \cite{chen2023}. Specifically, it accounts for generator output constraints on minimum and maximum output, minimum runtime constraints, maximum daily starts, and generator ramping constraints as discussed in \cite{vanackooij2018,yang2022,tahanan2015}. It further includes power flow balance and thermal constraints on transmission elements, such as lines and transformers. Unlike the linearized transmission constraints approaches discussed above, we use a continuous solver which enables us to adopt a classical continuous AC power flow formulation for transmission constraints \cite{bienstock2019}.

\cref{tab:parametersets} describes the input parameter sets of our UC optimization problem. The set of available generators $\mathcal{G}$ is comprised of the set of must-run generators $\mathcal{G}_m$ and the set of discretionary generators $\mathcal{G}_d$. So $G_m\cup G_d=G$. These sets change between iterations of the optimization problem, while the other sets stay fixed.

\begin{table}[htbp]
    \centering
    \caption{Input parameter sets for the optimization solver}
    \label{tab:parametersets}
    
    \renewcommand{\arraystretch}{1.1}
    \noindent
    \begin{tabularx}{0.5\textwidth}{@{} c c X @{}}
    \toprule
    \textbf{Notation} & \textbf{Indexing} &\multicolumn{1}{c}{\textbf{Meaning}} \\\midrule
    \(\mathcal{N}\) & & Set of buses \\\midrule
    \(\mathcal{G}\subseteq\mathcal{N}\) & & Set of must-run and discretionary generators \\\midrule
    \(\mathcal{G}_m\subseteq\mathcal{G}\) & & Set of must-run generators \\\midrule
    \(\mathcal{G}_d\subseteq\mathcal{G}\) & & Set of discretionary generators \\\midrule
    \(\mathcal{S}\) & & Set of contingency scenarios \\\midrule
    \(\mathcal{L}\) & & Set of lines \\\midrule
    \(\mathcal{L}_\delta\subseteq\mathcal{L}\) & \(\delta\in\mathcal{S}\) & Set of lines in contingency \(\delta\) \\\midrule
    \(N(i)\) & & Set of buses connected to bus \(i\) \\\midrule
    \(N_\delta(i)\) & \(\delta\in\mathcal{S}\) & Set of buses connected to bus \(i\) in contingency \(\delta\) \\\bottomrule
    \end{tabularx}
\end{table}

The decision variables of our problem are described in \cref{tab:decision variables}. The $u_i$s and $P_i$s are the variables the decision maker would need to commit units and decide which capacity they should produce at. The other variables are the result of the underlying physics, given these decisions. The corresponding variables with subscript \(\delta\) represent the values of the variables during contingency \(\delta\). \(P_{R,\delta,i}\) and \(Q_{R,\delta,i}\) are explicitly listed, since \(Q_{R,\delta,i}\) has no non-contingency counterpart, and \(P_{R,\delta,i}\)'s non-contingency counterpart is not a decision variable.

\begin{table}[htbp]
    \centering
    \caption{Decision Variables for the optimization solver}
    \label{tab:decision variables}
    
    \renewcommand{\arraystretch}{1.1}
    \noindent
    \begin{tabularx}{0.5\textwidth}{@{} c c X @{}}
    \toprule
    \textbf{Notation} & \textbf{Indexing} &\multicolumn{1}{c}{\textbf{Meaning}} \\\midrule
    \(u_i\) & \(i\in\mathcal{G}_d\) & Commitment status \\\midrule
    \(P_i\) & \(i\in\mathcal{G}\) & Real power production \\\midrule
    \(Q_i\) & \(i\in\mathcal{G}\) & Reactive power production \\\midrule
    \(V_i\) & \(i\in\mathcal{N}\) & Voltage magnitude \\\midrule
    \(\theta_i\) & \(i\in\mathcal{N}\) & Voltage phase angle (note: \(\theta_{ij} \doteq \theta_i - \theta_j, ij\in\mathcal{L}\)) \\\midrule
    \(P_{ij}\) & \(ij\in\mathcal{L}\) & Real power from bus \(i\) to bus \(j\) as measured at bus \(i\) \\\midrule
    \(Q_{ij}\) & \(ij\in\mathcal{L}\) & Reactive power from bus \(i\) to bus \(j\) as measured at bus \(i\) \\\midrule
    \(P_{R,\delta,i}\) & \(i\in\mathcal{N}, \delta\in\mathcal{S}\) & Real power from a generator that is ramping up at bus \(i\), during contingency \(\delta\)\\\midrule
    \(Q_{R,\delta,i}\) & \(i\in\mathcal{N}, \delta\in\mathcal{S}\) & Reactive power from a generator that is ramping up at bus \(i\), during contingency \(\delta\)\\
    \bottomrule
    \end{tabularx}
\end{table}

 Together with \cref{tab:parameters} which contains the optimization solver parameters, we can define our optimization problem:

 \begin{subequations}
\label{eq:UnitCommitment}
\begin{align}
\label{eq:Objective}
\min \; &
\sum_{i \in \mathcal{G}_m} \left(a_i P_i^2 + b_i P_i + c_i\right) 
\nonumber\\&+ \sum_{i \in \mathcal{G}_d} \left(a_i P_i^2 + b_i P_i + c_i\right) u_iu_{t-1,i}
\nonumber\\&+ \sum_{i \in \mathcal{G}_d}K_i(u_{t-1,i}(1-u_i)+(1-u_{t-1,i})u_i) \nonumber\\
& + \gamma\sum_{\delta\in\mathcal{S}}\sum_{i\in \mathcal{G}}u_i\left(a_i P_{\delta,i}^2 + b_i P_{\delta,i} + c_i\right)
\\[6pt]
\text{s.t.} \;
& \forall i\in\mathcal{G} : \max\!\left(P_{\min,i},\, P_{t-1,i} - r_{d,i}\right) \le P_i \nonumber\\&
\le \min\!\left(P_{\max,i},\, P_{t-1,i} + r_{u,i}\right)
\label{eq:ConsRamp}
\\&
\forall i\in\mathcal{G} : Q_{\min,i} \leq Q_i \leq Q_{\max, i}
\label{eq:ConsMVar}
\\&
\forall i\in\mathcal{G} : u_i\in \{0,1\}
\label{eq:ConsBools}
\\&
\forall i \in \mathcal{N} : \;(P_{i}u_i + P_{\min,i}(1-u_i))u_{t-1,i} + P_{R,i} - P_{Di} \nonumber\\& =G_{\text{sh},i} + \sum_{j \in N(i)} P_{ij} \nonumber\\& \doteq
\sum_{j \in N(i)\cup\{i\}}  V_{i}V_{j}
\left( G_{ij}\cos\theta_{ij} + B_{ij}\sin\theta_{ij} \right)
\label{eq:ConsPflow}
\\&
\forall i \in \mathcal{N} : \;Q_{i}u_iu_{t-1,i} - Q_{Di}=-B_{\text{sh},i} + \sum_{j \in N(i)} Q_{ij} \nonumber\\& \doteq 
\sum_{j \in N(i)\cup\{i\}}V_{i}  V_{j}
\left( G_{ij}\sin\theta_{ij} - B_{ij}\cos\theta_{ij} \right)
\label{eq:ConsQflow}
\\&
\forall ij\in\mathcal{L} : P_{ij}^2 + Q_{ij}^2 \le S_{\max,ij}^2
\label{eq:ConsThermal}
\\&
\forall i\in\mathcal{N} : V_{\min,i} \le V_{i} \le V_{\max,i}
\label{eq:ConsVoltage}
\\&
\forall ij\in\mathcal{L} : \theta_{ij} = \theta_{i} - \theta_{j}
\label{eq:ConsAngle}
\\&
\forall i\in\mathcal{G} : P_{\min,i} \leq P_{\delta,i},P_{R,\delta,i} \leq P_{\max,i}
\label{eq:ConsPDelta}
\\&
\forall i\in\mathcal{G} : Q_{\min,i} \leq Q_{\delta,i},Q_{R,\delta,i} \leq Q_{\max,i}
\label{eq:ConsQDelta}
\\&
\forall i \in \mathcal{N}, \;\delta\in \mathcal{S}: \; P_{\delta,i}u_i + P_{R,\delta,i} - P_{\delta,Di}\nonumber = \nonumber\\&G_{\text{sh},i}+\sum_{j \in N_\delta(i)} P_{\delta,ij} \nonumber\\& \doteq
\sum_{j \in N_\delta(i)\cup\{i\}}  V_{\delta,i}V_{\delta,j}
\left( G_{ij}\cos\theta_{\delta,ij} + B_{ij}\sin\theta_{\delta,ij} \right)
\label{eq:ConsPflowdel}
\\&
\forall i \in \mathcal{N}, \delta\in \mathcal{S}: \; Q_{\delta,i}u_i + Q_{R,\delta,i} - Q_{\delta,Di}\nonumber\\& = -B_{\text{sh},i}+\sum_{j \in N_\delta(i)} Q_{\delta,ij} \nonumber\\& \doteq
\sum_{j \in N_\delta(i)\cup\{i\}}V_{\delta,i}  V_{\delta,j}
\left( G_{ij}\sin\theta_{\delta,ij} - B_{ij}\cos\theta_{\delta,ij} \right)
\label{eq:ConsQflowdel}
\\&
\forall \delta\in\mathcal{S}, ij\in\mathcal{L}_\delta : P_{\delta,ij}^2 + Q_{\delta,ij}^2 \le S_{\max,ij}^2
\label{eq:ConsThermaldel}
\\&
\forall \delta\in\mathcal{S}, i\in\mathcal{N} : V_{\min,i} \le V_{\delta,i} \le V_{\max,i}
\label{eq:ConsVoltagedel}
\\&
\forall \delta\in\mathcal{S}, ij\in\mathcal{L}_\delta : \theta_{\delta,ij} = \theta_{\delta,i} - \theta_{\delta,j}
\label{eq:ConsAngledel}
\end{align}
\end{subequations}

\begin{table}[htbp]
    \centering
    \caption{Input Parameters for the optimization solver}
    \label{tab:parameters}
    
    \renewcommand{\arraystretch}{1.1}
    \noindent
    \begin{tabularx}{0.5\textwidth}{@{} c c X @{}}
    \toprule
    \textbf{Notation} & \textbf{Indexing} &\multicolumn{1}{c}{\textbf{Meaning}} \\
    \(\gamma\) & & Bias factor to select units that will be more efficient at meeting contingency load \\\midrule
    \(a_i,b_i,c_i\) & \(i\in\mathcal{G}\) & Quadratic cost coefficients, usually resulting from fuel costs or strategic bidding behavior \\\midrule
    \(K_i\) & \(i\in\mathcal{G}\) & Commitment change cost - The true startup and shutdown costs are rarely the same. Averaging them prevents market distortion and removes artificial incentives to turn on or off.\\\midrule
    \(u_{t-1},i\) & \(i\in\mathcal{G}\) & Previous state (1 \(\leftrightarrow\) on, 0 \(\leftrightarrow\) off) \\\midrule
    \(P_{t-1},i\) & \(i\in\mathcal{G}\) & Previous real power output \\\midrule
    \(r_{u,i}, r_{d,i}\) & \(i\in\mathcal{G}\) & Real power ramp-up and ramp-down limits per timestep \\\midrule
    \(P_{\min,i},P_{\max,i}\) & \(i\in\mathcal{G}\) & Minimum and maximum real power output \\\midrule
    \(Q_{\min,i},Q_{\max,i}\) & \(i\in\mathcal{G}\) & Minimum and maximum reactive power output \\\midrule
    \(P_{R,i}\) & \(i\in\mathcal{N}\) & Real power from ramping generators at bus \(i\)\\\midrule
    \(P_{Di}\) & \(i\in\mathcal{N}\) & Real power load \\\midrule
    \(Q_{Di}\) & \(i\in\mathcal{N}\) & Reactive power load \\\midrule
    \(P_{\delta,Di}\)& \(\delta\in\mathcal{S},i\in\mathcal{N}\) & Real power load in contingency \(\delta\) \\\midrule
    \(Q_{\delta,Di}\) & \(\delta\in\mathcal{S},i\in\mathcal{N}\) & Reactive power load in contingency \(\delta\) \\\midrule
    \(V_{\min,i},V_{\max,i}\) & \(i\in\mathcal{N}\) & Minimum and maximum voltage magnitude \\\midrule
    \(S_{\max,ij}\) & \(ij\in\mathcal{L}\) & Line thermal limit \\\midrule
    \(G_{ij}, B_{ij}\) & \(ij\in\mathcal{L}\) & Real and imaginary admittance values \\\midrule
    \(G_{\text{sh},i} = \sum_{j:ij\in\mathcal{L}}G_{ij}\) & \(i\in\mathcal{N}\) & Real power shunt value \\\midrule
    \(B_{\text{sh},i} = \sum_{j:ij\in\mathcal{L}}B_{ij}\) & \(i\in\mathcal{N}\) & Reactive power shunt value \\
    \bottomrule
    \end{tabularx}
\end{table}

In our implementation, we set \(\gamma = 1\) in \cref{eq:Objective}, and consider two contingency scenarios, which are the minimum and maximum load in a 72-hour prediction window (\(|\mathcal{S}|=2\)). In line with the area interchange logic used by large independent system operators, we subdivide the system into \textit{areas}, while maintaining MVA constraints on each transmission element so that the model remains nodal. Each resulting area has relatively high internal connectivity and weaker connectivity to the rest of the system. 

We account for loss of load and downward deviations from the load forecast by reducing minimum load in each area \(A\) by \(2\% + \max_{i\in A}\{P_{Di}\}/\sum_{i\in A}P_{Di}\). 
We account for upwards deviations from load forecast by increasing maximum load in each area \(A\) by \(5\% + \max_{i\in A\cap\mathcal{G}}\{P_{\max,i}\}/\sum_{i\in A\cap\mathcal{G}} P_{\max,i}\). For the extreme peak loads, we also test the loss of a transmission line with a base voltage of 345 kV (inter-area transmission corridors). Taking the maximum and minimum projected loads utilizes the cost curves' monotonicity. If a set of units can supply these edge cases, it can supply any load between them.

This choice is consistent with results from \cite{regev2026economicUC} which show that most of the gains of incorporating future costs effects into unit commitment come from the projected minimum and maximum loads.

Any bus that has no generator, or has a generator which cannot run, has its power production and commitment status variables constrained to 0. Mathematically, in \cref{eq:ConsPflow,eq:ConsQflow,eq:ConsPflowdel,eq:ConsQflowdel} \(P_i=P_{\delta,i}=Q_i=Q_{\delta,i}=u_i=0\) for \(i\in\mathcal{N}\backslash\mathcal{G}\).

\section{Simulation design}
\label{sec:design}

\subsection{Transmission constraints}

\begin{figure}[htbp]
\centerline{\includegraphics[width=\columnwidth]{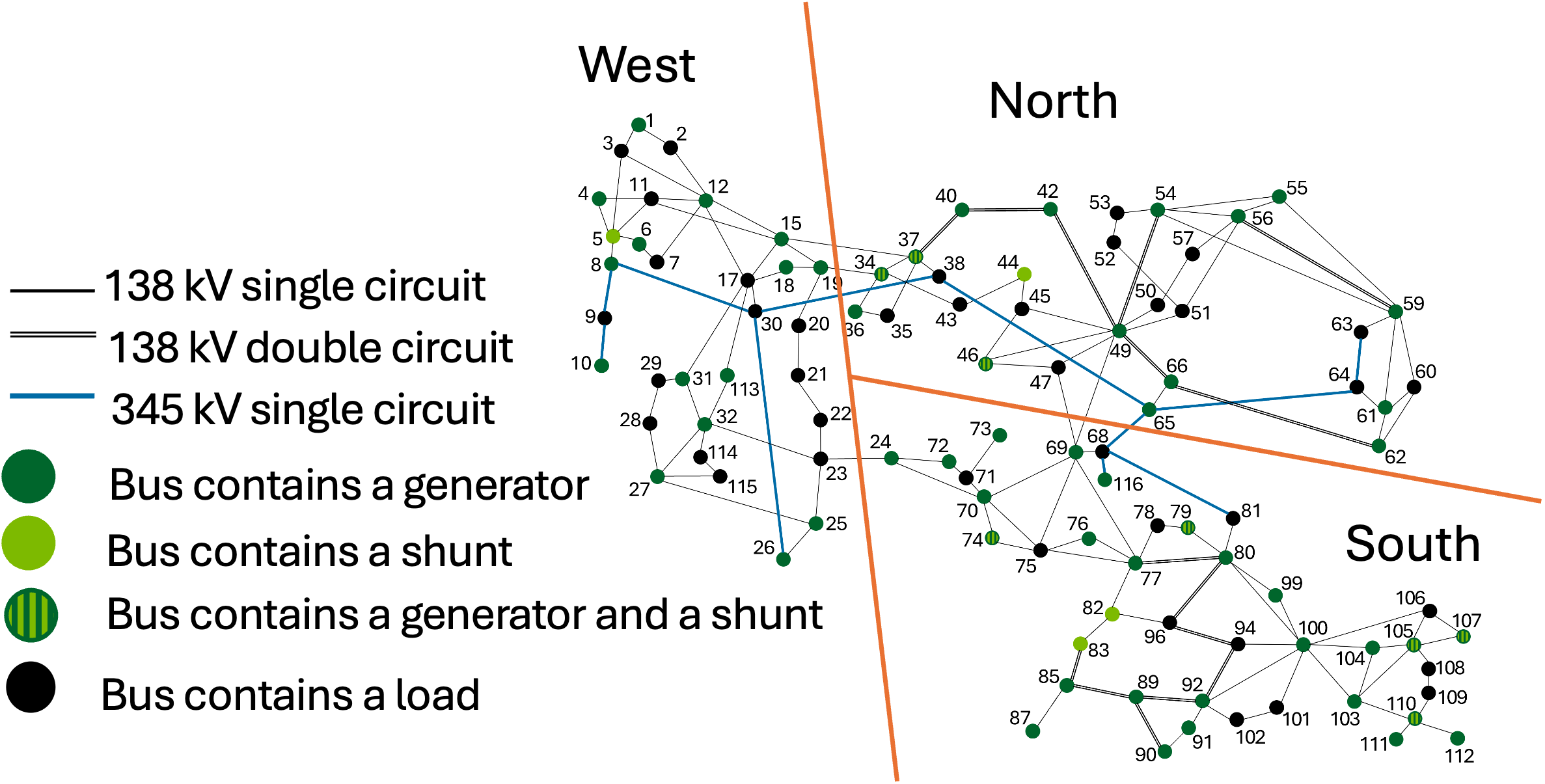}}
\caption{Modified IEEE 118-bus grid \cite{tsybina2026agenticartificialintelligencepower}.}
\label{fig:adjusted_topology}
\end{figure}

To introduce transmission constraints, we use an adjusted version of the IEEE 118 bus system \cite{ieee118bus}, modified in two major ways \cite{tsybina2026agenticartificialintelligencepower}. First, we reduce the number of buses from 118 to 100 for convenient increments. Second, we apply transmission and generation limits from \cite{PNNL_grid}, which are strict enough to create a congestible grid. The resulting system represents a three-area model with two voltage classes (138 kV and 345 kV), with comparatively easy intra-area exchange and with inter-area exchange that is constrained by MVA flows through multi-line transmission corridors. This configuration of the system allows UC to select the lowest cost configuration of generating capacity assuming nodal level of power flow detail, while allocating reserves per area. \cref{fig:adjusted_topology} illustrates the adjusted grid configuration. 

\begin{figure}[htbp]
\centerline{\includegraphics[width=\columnwidth]{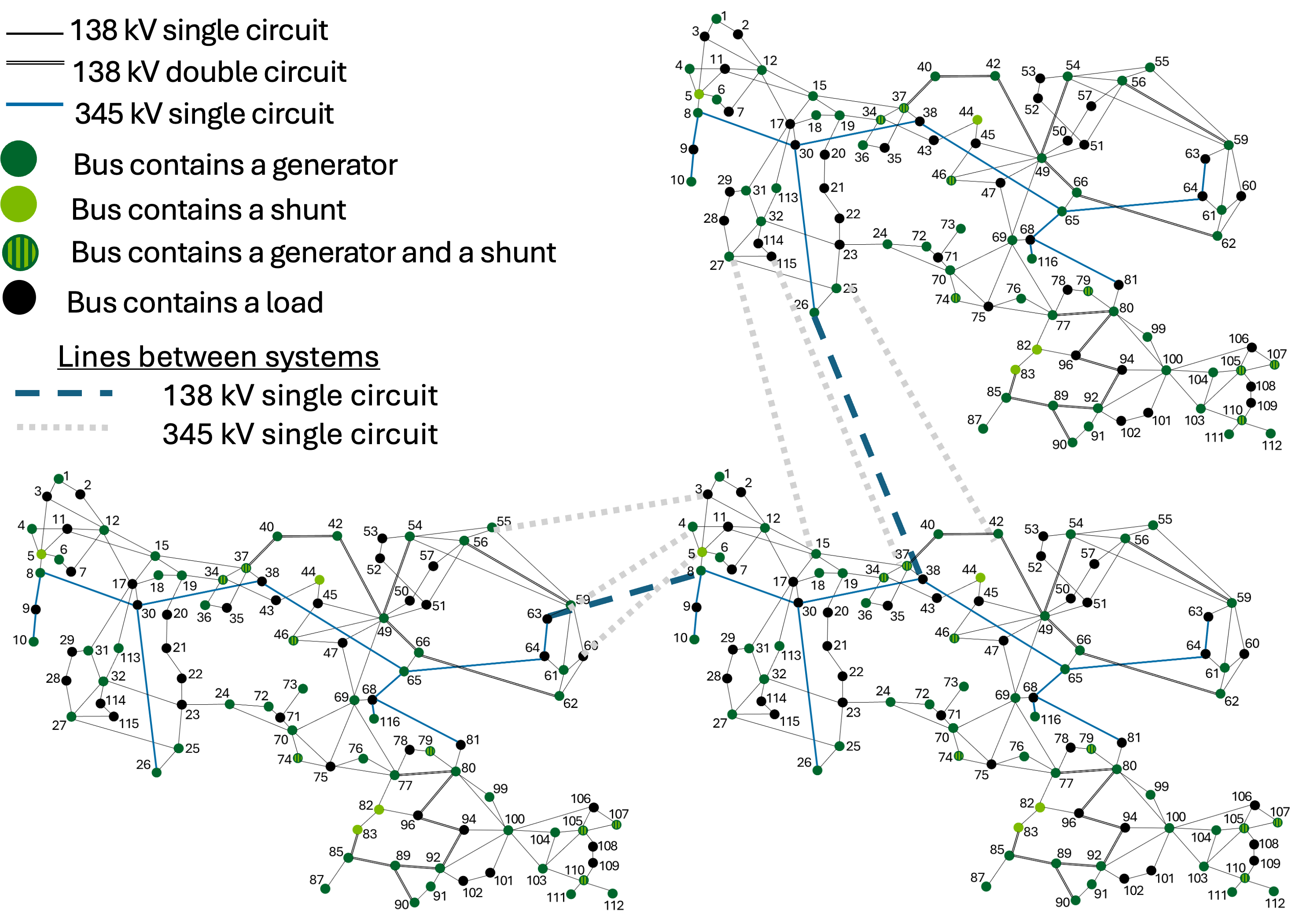}}
\caption{Horizontal and vertical extension of \cref{fig:adjusted_topology}, the modified IEEE 118-bus grid.}
\label{fig:extended_topology}
\end{figure}

The system is designed to be modular and extends horizontally or vertically in an arbitrary number of steps. The vertical extension uses the connection of buses 15-27, 37-115, 42-25, 38-26 between grid copies. The horizontal extension uses the connection of buses 3-55, 4-59, 5-60, 8-63. The grid extension creates additional lines with the circuit characteristics similar to those assigned to existing lines in the system, for the respective voltage class. We extend the system to create total system sizes ranging from 100 to 20,000 buses (which is the size of the largest RTO in the US, PJM). Each extension of the grid is added in a random direction and perturbed within \(\pm1\%\) to avoid symmetry in load and generation patterns which could appear if the grid were extended through \(n\) identical copies.  \cref{fig:extended_topology} shows an example with a vertical and horizontal extension of \cref{fig:adjusted_topology}.

% We test each contingency per-area (West, North, South) within each copy of the grid. 

\subsection{Generation assumptions}

We test the RRUC algorithm using the generator characteristics available from \cite{PNNL_grid}, augmented with ramp and runtime constraints. The single 100 bus model includes 54 generators, each in a separate bus.  

The original IEEE model and operating characteristics from \cite{PNNL_grid} do not include runtime constraints and ramp up and ramp down information, which are essential for UC. We used the same logic as in~\cite{regev2026economicUC} to assign operating constraints by fuel type and unit size. \cref{tab:datarun} shows the additional assumptions about equipment characteristics. The information on runtime constraints and ramp characteristics is available from a number of US and EU sources \cite{diw2013, GE}. Warm start costs are available from PJM~\cite{pjm} and validated against the EU system~\cite{diw2013}. 
Ramp characteristics are provided in \cite{PNNL_grid}, but we opted for manufacturer specifications to inform ramp characteristics in order to preserve the consistency of mechanical generator data for the same unit. These numbers  are likely an upper bound of ramp rates. Evidence from real systems suggests that power plants do not always deliver the nominal ramp rates specified for equipment~\cite{India}.

\begin{table}[htb]
\caption{\label{tab:datarun} Runtime, Ramp Assumptions, and Data Sources}
\begin{tabular}{lrrr} \toprule
   Parameter & Source & Coal & Gas
\\ \midrule
Min runtime [min] & ~\cite{pjm} & 240-1440 & 0-1440
\\ Max daily starts & ~\cite{pjm} & 1-3 & 1-24
%\\ Hot start [min] & ~\cite{diw2013,GE} & 60-240 & 25-120
\\ Warm start [min] & ~\cite{diw2013} & 120-480 & 60-240
%\\ Cold start [min] & ~\cite{diw2013} & 360-720 & 120-300
% \\ Relative cost of start, warm/hot & ~\cite{pjm,diw2013} & 1.19-1.42 & 1.12-1.58
% \\ Relative cost of start, cold/hot & ~\cite{pjm,diw2013} & 1.74-1.93 & 1.35-2.25 \\
\\Ramp up rate [\%/min] & ~\cite{diw2013,GE} & 1-6 & 2-12
\\ Ramp down rate [\%/min] & ~\cite{diw2013} & 5 & 15
\\ \bottomrule
\end{tabular}
\end{table}

% \textcolor{red}{ADD HYDRO ASSUMPTIONS.} 
The original data  does not have information about generator cost curves. We use prices in individual generator bids for June 2025~\cite{pjm} and match it to \cite{PNNL_grid} based on \(P_{\min}\), \(P_{\max}\). First, we match each generator to the bid data based on the generator's fuel type. Then, we use the \(P_{\min}\) and \(P_{\max}\) to inform the closest match for operating characteristics. The bid data from \cite{pjm} does not offer specific fuel types for generators. We inferred fuel type of \cite{pjm} units from the steps of supply curves and the runtime characteristics, which included natural gas and coal generators. 

\begin{figure}[htbp]
\centerline{\includegraphics[width=\columnwidth]{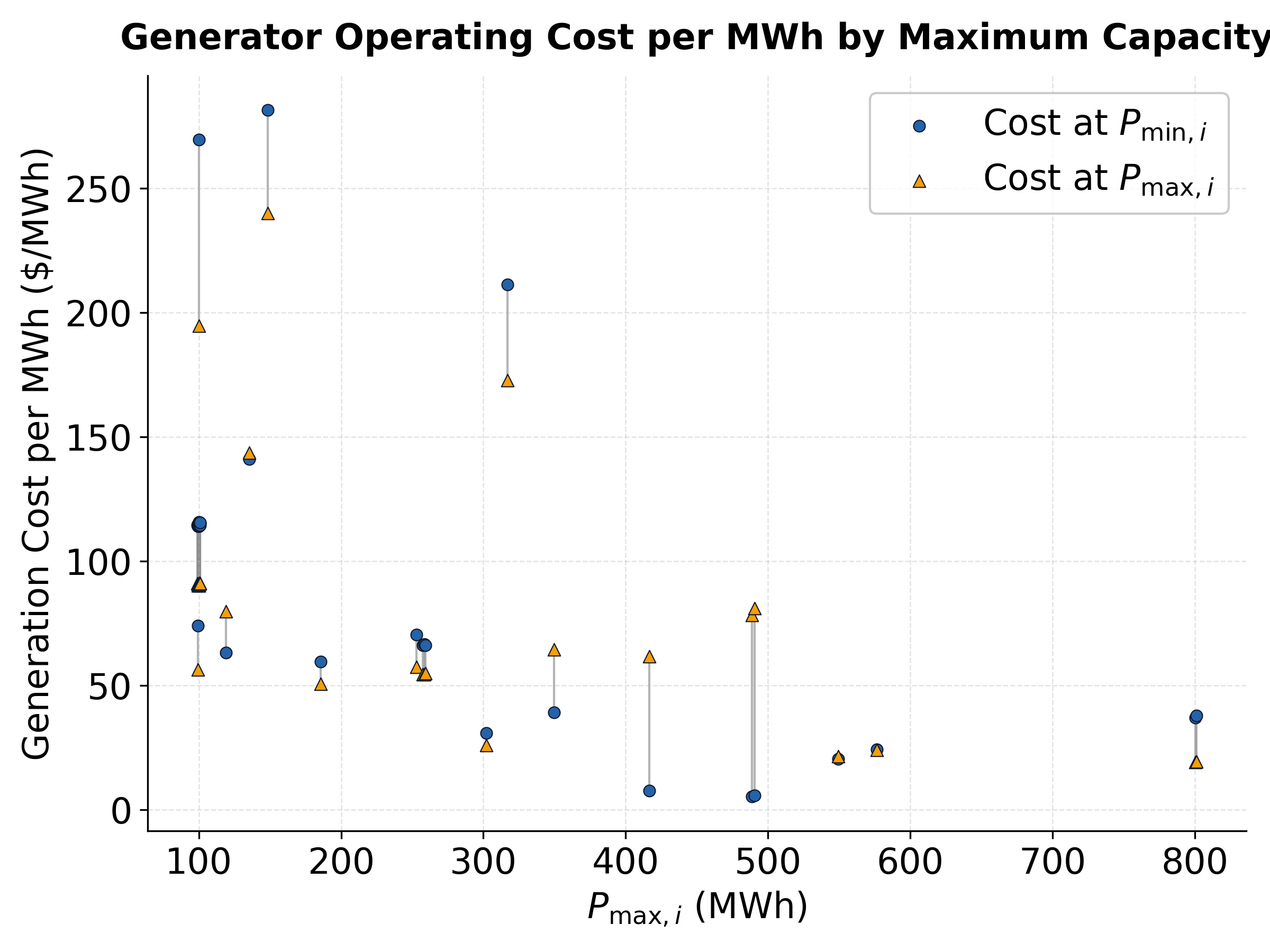}}
\caption{Cost per Megawatt-hour at maximum and minimum capacity based on generator size (\(P_{\max,i}\)).}
\label{fig:gen_data}
\end{figure}

Generators 24, 85, 87, 89, 90, 91, and 92 did not have fuel associated with them in \cite{PNNL_grid}. For those generators, we used \(P_{\min}\) to assign fuel type, which was consistent with all but 3 of the units with an associated fuel type. If \(P_{\min}<25\) MW, the unit was assigned as hydro. If \(P_{\min}\in[25, 80)\) MW, it was assigned natural gas.
%Units with combined cycles were assumed to be gas turbines. 
Generators in which \(P_{\min}\geq 80\) MW were assigned coal. Once all generators were matched with bid shapes from \cite{pjm}, we used a postprocessing step to convert piecewise linear curves from \cite{pjm} into twice differentiable cost curves for natural gas and coal generators. 
\cref{fig:gen_data} presents the distribution of generator cost at \(P_{\min,i}\), cost at \(P_{\max,i}\), and unit size, excluding the 19 hydro generators.

We assumed hydro generators were large reservoir units, with enough rainfall to operate consistently. We chose $a,b$, the quadratic and linear cost coefficients, for hydro generators randomly in the range  \([10^{-7},10^{-6}] \) and set $c,K$, the constant term and on/off cost, to 0. This models the cost of hydro power production as essentially free. The resulting dataset had coherent \(P_{\min}\), \(P_{\max}\), runtime and ramp constraints, and cost curves.

\subsection{Load assumptions}

We obtained a load profile from historical PJM records~\cite{pjm}. The PJM load profile is applied to IEEE 188 bus system load and scaled by bus, assuming constant power factor. The historical high load was recorded in PJM on June 23, 2025 at 5 pm, with total load of 160.2 GW~\cite{pjm}. We study one week around this peak. The resulting runtime constrained problem has 672 time periods at 15 minute intervals (7 days). 

The simulation starts four days before the peak to ensure that the algorithm has enough room to adjust. The unit commitment output for the first day is used as a warm-up for the system, to avoid cold starting the grid, the results for this day are discarded. Simulation ends 72 hours before the end of the study period, so that the last solve provides a 72 hour committed generation set for the remaining load forecast. This means that, out of the original 7-day window, the results displayed in \cref{sec:results} use days 2, 3, and 4 (periods 97 through 384).

\section{Implementation and Algorithm}
\label{sec:algorithm}

We implemented our algorithm in Julia~\cite{Bezanson2017Julia} using the MadNLP nonlinear programming solver~\cite{shin2021graph, shin2024accelerating}, which built upon~\cite{Regev2022HyKKT, Swirydowicz2025GPU}. We use the JuMP modeling language~\cite{Dunning2017JuMP}, and the MA57 linear solver from HSL~\cite{HSL}. We solve all optimization problems to $10^{-3}$ relative tolerance. \footnote{This is a common practice for line constrained problems~\cite{saadat2010power}. We do this also to give Gurobi an opportunity to arrive at a solution. Our solver is able to solve problems to $10^{-6}$ tolerance.}  We note this is not a limitation of our solver and only decreased\cref{alg:rruc} shows pseudocode for the algorithm.

\begin{algorithm}[htbp]
\caption{Topologically Constrained RRUC}
\label{alg:rruc}
\begin{algorithmic}[1]
\Function{RoundSolution}{\mbox{$\{y_i^*\}$}}
    \State \textbf{assert} \(|\{y_i^*\}| = |\mathcal{G}_d|\)
    \State Re-index such that \(y_1^*\geq y_2^*\geq\ldots\geq y_n^*\)
    \State Attempt to solve the optimal flow problem induced by rounding up all \(y_i^* > 0\)
    \State\textbf{if} solving succeeds \textbf{return} solution
    \For{\(\tau \gets \max\{i:y_i^*=1\}\) \textbf{to} \(|\mathcal{G}_d|\)}
        \State Attempt to solve the optimal flow problem induced by rounding up \(\{y_{i\leq\tau}^*\}\)
        \State\textbf{if} solving succeeds \textbf{return} solution
    \EndFor
\EndFunction
\Statex
\For{\textnormal{period \(\gets 1\)} \textbf{to} \textnormal{num\_periods}}
    \State Construct relaxed instance \(\mathcal{I}\) of \cref{eq:UnitCommitment} with
    \Statex \(u_i\in\{0,1\}\rightarrow y_i\in[0,1]\)
    \State Solve \(\mathcal{I}\)
    \State \(y_i^*\gets \) resulting assignments to \(y_i\)
    \State Save \Call{RoundSolution}{\mbox{$\{y_i^*\}$}}
    \State Update generator states, ramping power, and \(\mathcal{G}_d\), \(\mathcal{G}_m\)
\EndFor
\end{algorithmic}
\end{algorithm}

Our implementation optionally uses JACC~\cite{JACC} for parallelization on lines 6-9, though this did not accelerate the solver on our system with the current implementation. \cref{alg:rruc} is parallelizable in theory, a more specialized implementation may lead to parallel speedups. We save this for future work.

\section{Results}
\label{sec:results}

\subsection{Computational efficiency}

\begin{figure}[htbp]
    \centerline{\includegraphics[width=\columnwidth]{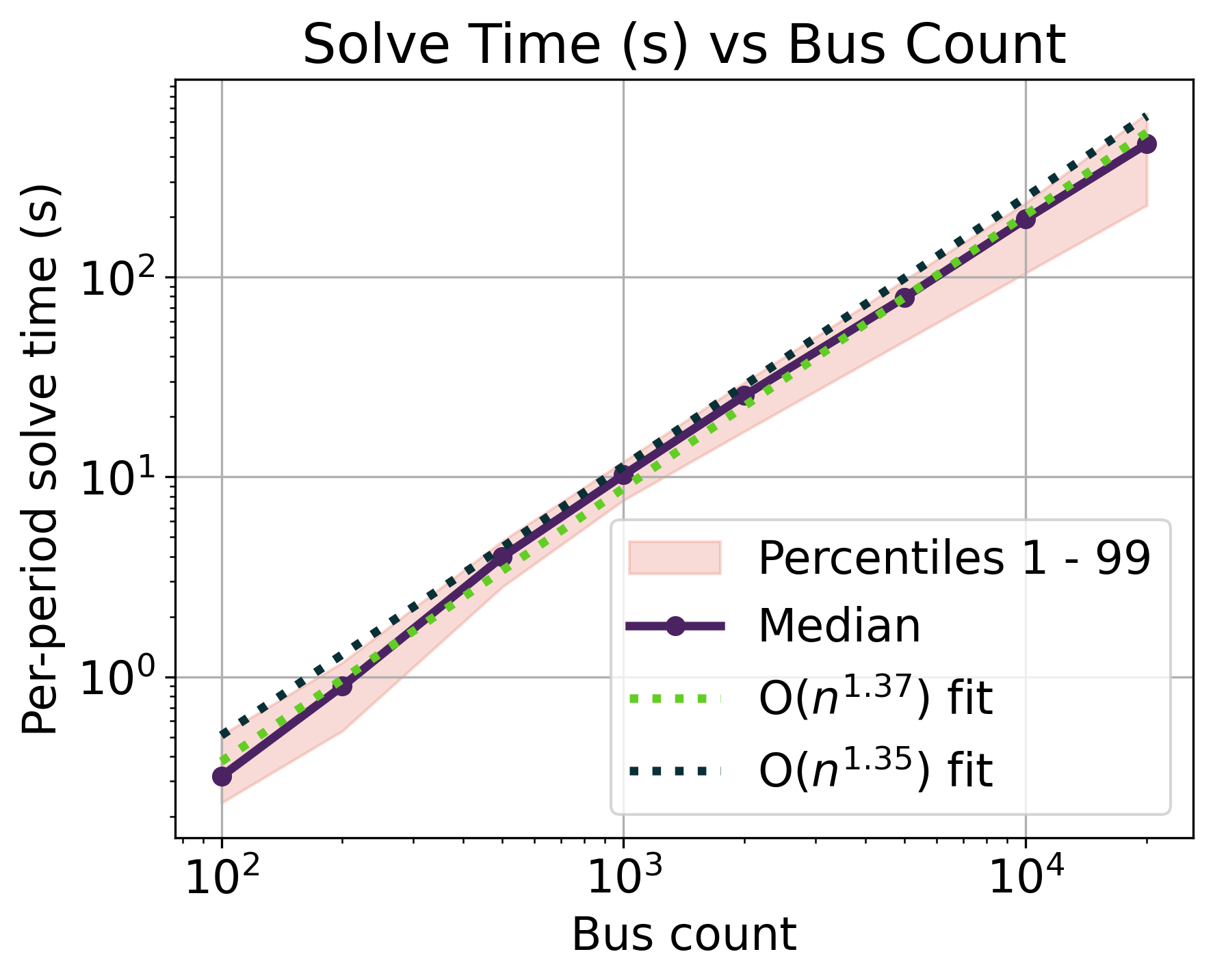}}
    \caption{Experiment results for per-period solve time. RRUC's mean solve time scales as \(\mathcal{O}\left(n^{\text{1.35}}\right)\) in the numbers of buses. The relatively large deviation to the 1 percentile solve time is due to a few periods where the problem is more constrained because of mechanical limits on the generators.}
    \label{fig:timing_results}
\end{figure}

\cref{fig:timing_results} summarizes the runtime scaling results of our  computational experiments. It shows that the median run time per-period for our implementation scales as \(\mathcal{O}\left(n^{\text{1.35}}\right)\), where \(n\) is the number of buses in the system. The 99\textsuperscript{th} percentile solve time scales as \(\mathcal{O}\left(n^{1.37}\right)\). These deviations are caused by periods that have fewer (more) variables and constraints, because more (fewer) of their units are fixed. 

Gurobi~\cite{Gurobi} was not able to find a feasible solution to the problem within 15 minutes to the 100-bus system, while RRUC takes less than 1 second. This matches results from \cite{regev2026economicUC}, which show that Gurobi requires time exponential in the number of constraints. This constitutes a significant methodological advancement over \cite{regev2026economicUC}, where in the absence of transmission constraints Gurobi was competitive. Including transmission MVA limits and power flow balance equations increases the number of constraints linearly in the number of buses. Even the smallest 100-bus system has hundreds of these constraints. This explains why Gurobi is not an efficient solver for our transmission constrained problem in \cref{eq:UnitCommitment}. 

\begin{figure}[htbp]
    \centerline{\includegraphics[width=\columnwidth]{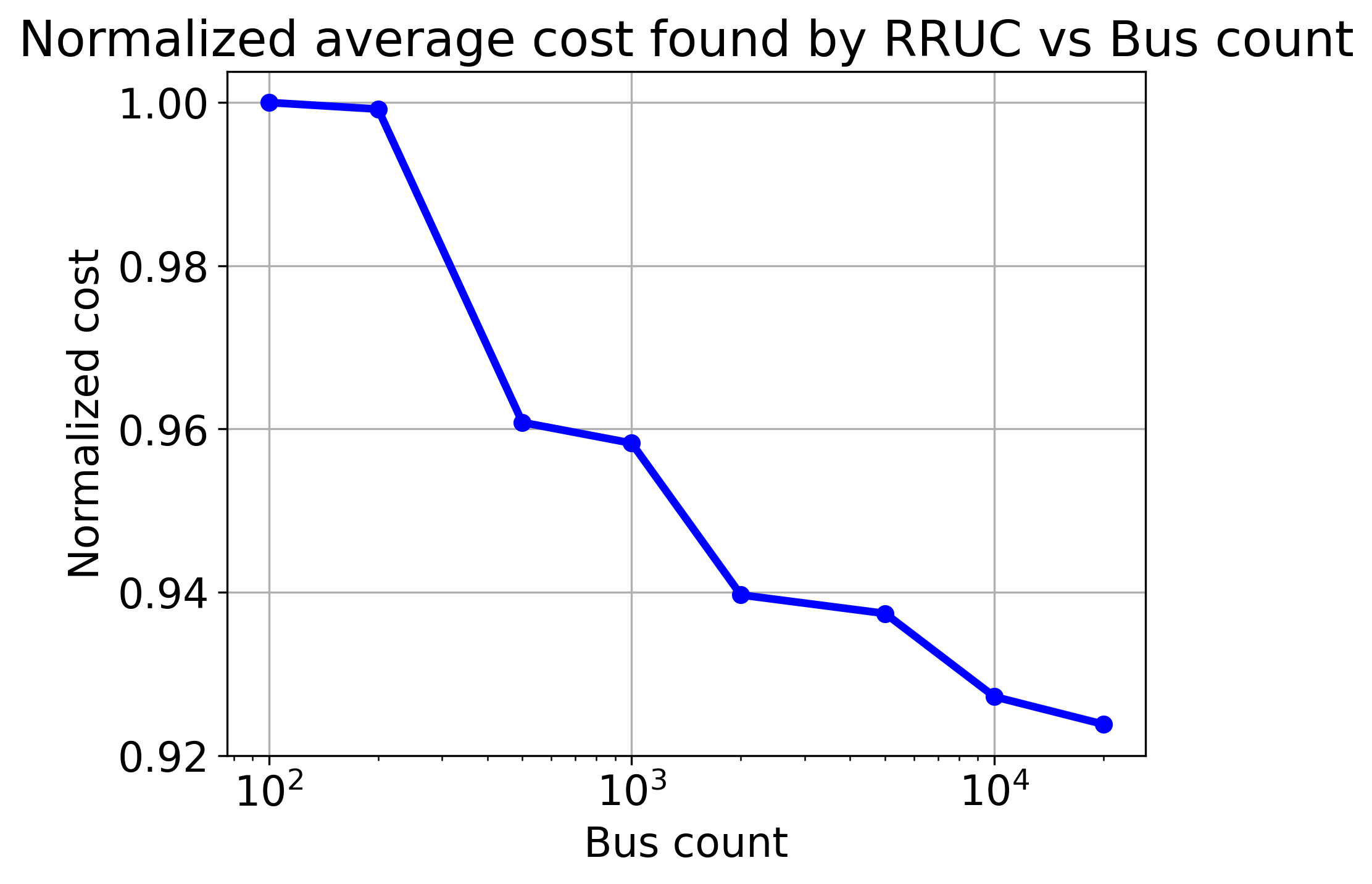}}
    \caption{Results for mean per-period operation cost, normalized to 1 for 100 buses. Optimizing the 20,000 bus grid simultaneously led to an almost 8\% cost reduction versus optimizing each 100 bus section separately. }
    \label{fig:cost_normalized}
\end{figure}

\cref{fig:cost_normalized} shows the average cost per time period, normalized by the number of buses, such that 100 buses has a cost of 1. The normalized cost steadily decreases as the problem grows, showing that optimizing over larger number of generators can greatly reduce cost. The savings are larger than the difference in production costs between the copies, indicating the redispatch effect. We expect these savings to grow when optimizing over larger systems with more heterogeneous generation and load than tested here. 

\begin{figure}[htbp]
    \centerline{\includegraphics[width=\columnwidth]{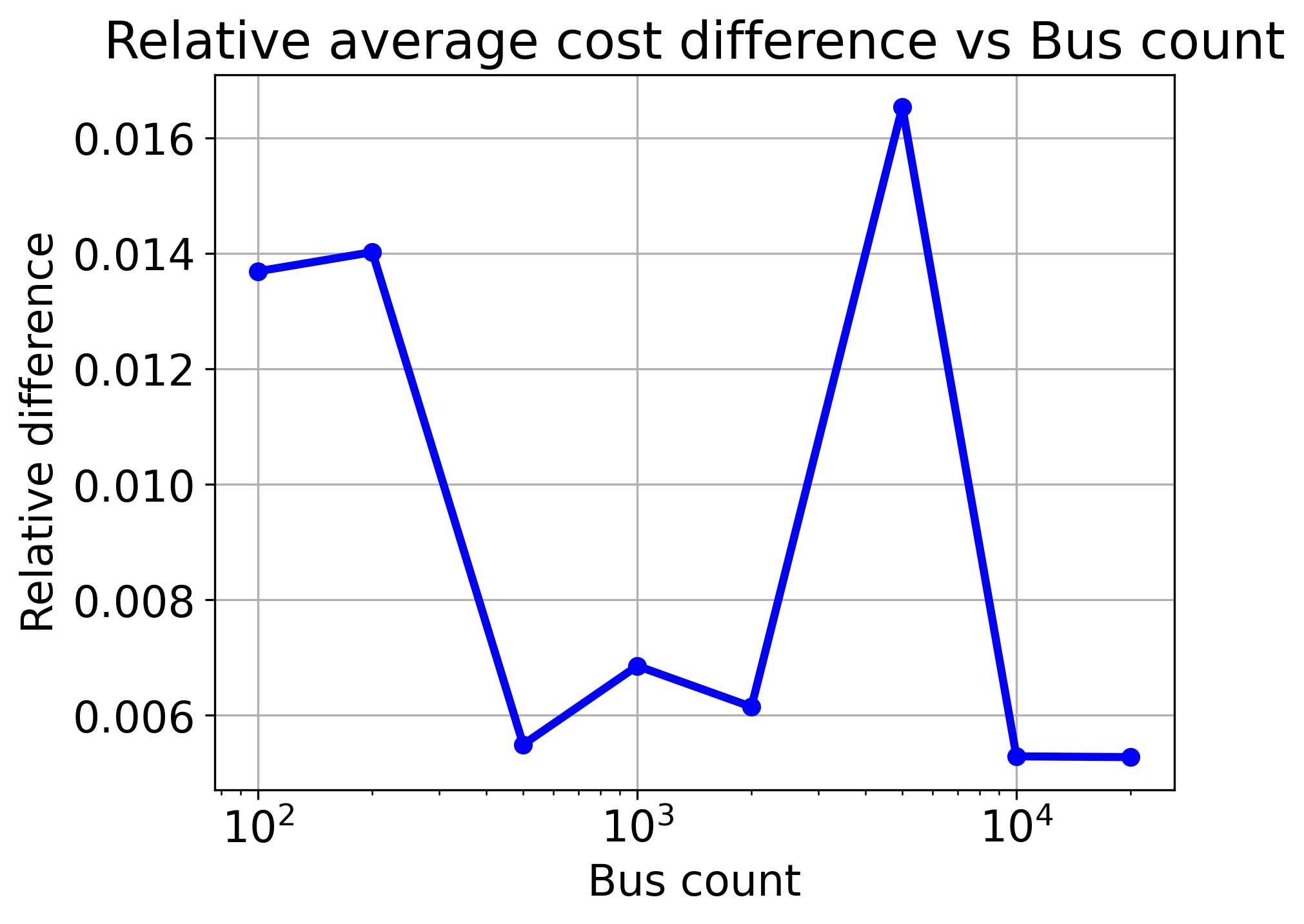}}
    \caption{$\mathbf{R}_d=(\mathbf{C}_C-\mathbf{C}_R)/\mathbf{C}_C$ with $\mathbf{C}_R$ the average production cost of the solution found by RRUC   and $\mathbf{C}_C$ that of the (perhaps infeasible) continuous relaxation found on Line 13 of \cref{alg:rruc}. $\mathbf{C}_R$ is slightly better, likely an artifact of the fact that we are only indirectly minimizing $\mathbf{C}_R$ and $\mathbf{C}_C$. We take this to mean that RRUC finds solutions that are numerically indistinguishable from the continuous relaxation.}
    \label{fig:cost_relative}
\end{figure}

To find the average cost of a solution, we average the objective \cref{eq:Objective} without the final term over all time periods. The last term in \cref{eq:Objective} is there to help lower production costs of future time periods, but is not a cost incurred in practice. Consistent with the results of \cite{regev2026economicUC}, our tests found that including this term reduced average production costs.
Let the average production cost of the solution found by RRUC be $\mathbf{C}_R$ and that of the continuous relaxation be $\mathbf{C}_C$ found on Line 13 of \cref{alg:rruc}. \cref{fig:cost_relative} shows the relative  difference between them $\mathbf{R}_d=(\mathbf{C}_C-\mathbf{C}_R)/\mathbf{C}_C$. 

$\mathbf{C}_C$ consistently larger than $\mathbf{C}_R$ is an artifact of the objective \cref{eq:Objective} over the current time period being a proxy for the real target \cref{eq:Objective} without the last term, averaged over all periods. The last term in \cref{eq:Objective} is typically larger than the other terms combined. This means that a relative error of $10^{-3}$ in \cref{eq:Objective} could be larger in relative terms without the last term. We take this to mean that the solutions are just as good for practical purposes, even though RRUC's average in \cref{fig:cost_relative} is lower than that of the (perhaps infeasible) continuous solution on Line 13 of \cref{alg:rruc}.

% If we assume that the optimal operation cost is proportional to \(n\), this then implies that our method has at most a small constant multiplicative error. That is, if \(\alg\) and \(\opt\) are the costs of our solution and an optimal solution, respectively, we have
% \begin{align*}
%     \alg = \theta(n),\ \opt = \theta(n) \Rightarrow \alg = \theta(\opt)
% \end{align*}

\subsection{Engineering reliability}

\begin{figure*}[!htbp]
    \centerline{\includegraphics[width=\textwidth]{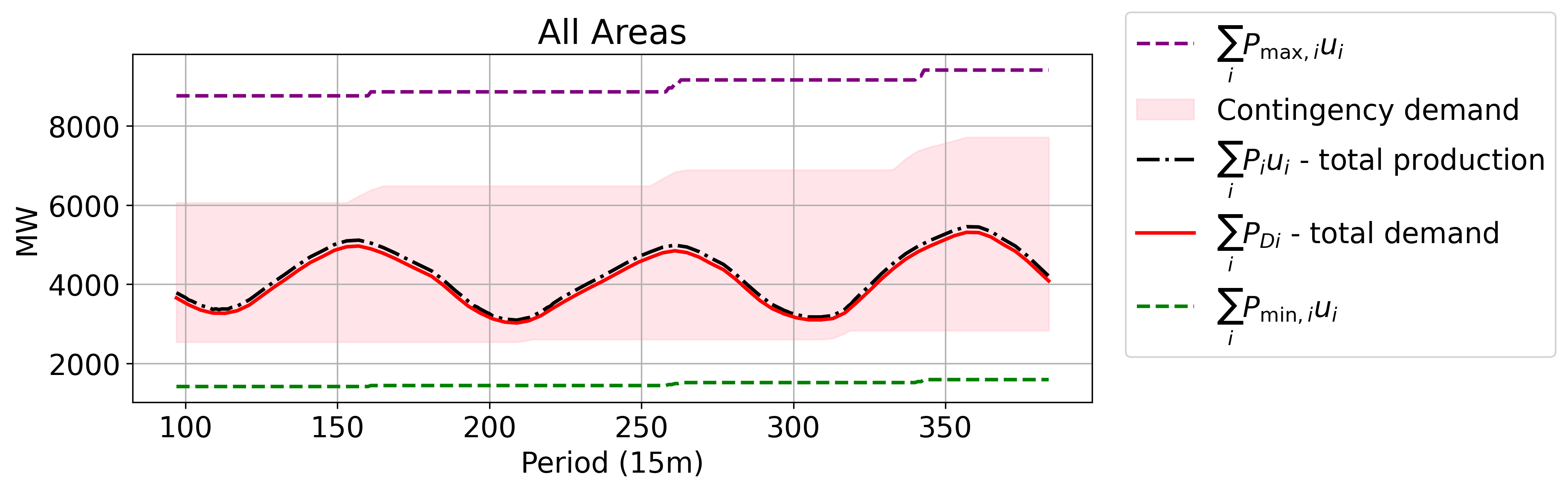}}
    \caption{Real power capacity, production, and load over all areas in the 100-bus system. Production tracks load. The selected generators can produce the minimum and maximum required load by the contingencies. The change in maximum and minimum production abilities (purple and green) occurs a few periods after new maximum or minimum contingencies enter the time window. The difference corresponds to the ramp up times. RRUC chooses to turn on more generators to more efficiently meet growing load in the future, even though it is not required to meet the contingency.}
    \label{fig:all areas}
\end{figure*}

\begin{figure*}[!htbp]
    \centerline{\includegraphics[width=\textwidth]{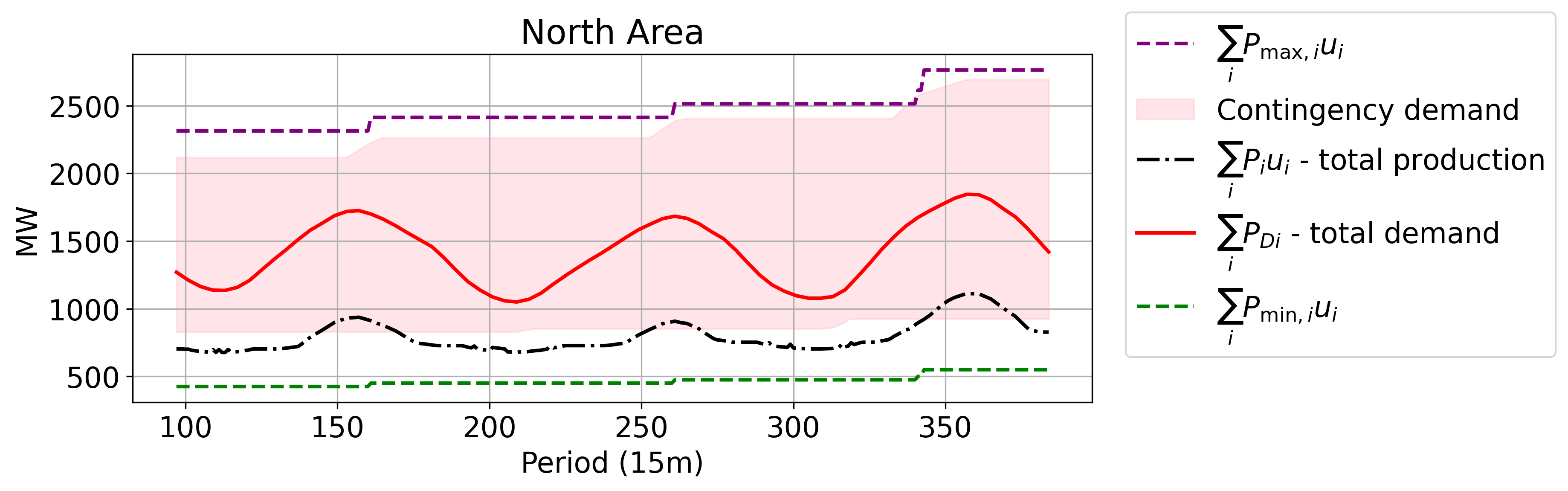}}
    \caption{Real power capacity, production, and load in the North area of the 100-bus system. The North area imports power, including to meet its contingencies, because it is allowed to do so by the constraints (which only take down some of the lines).}
    \label{fig:north}
\end{figure*}

\begin{figure*}[!htbp]
    \centerline{\includegraphics[width=\textwidth]{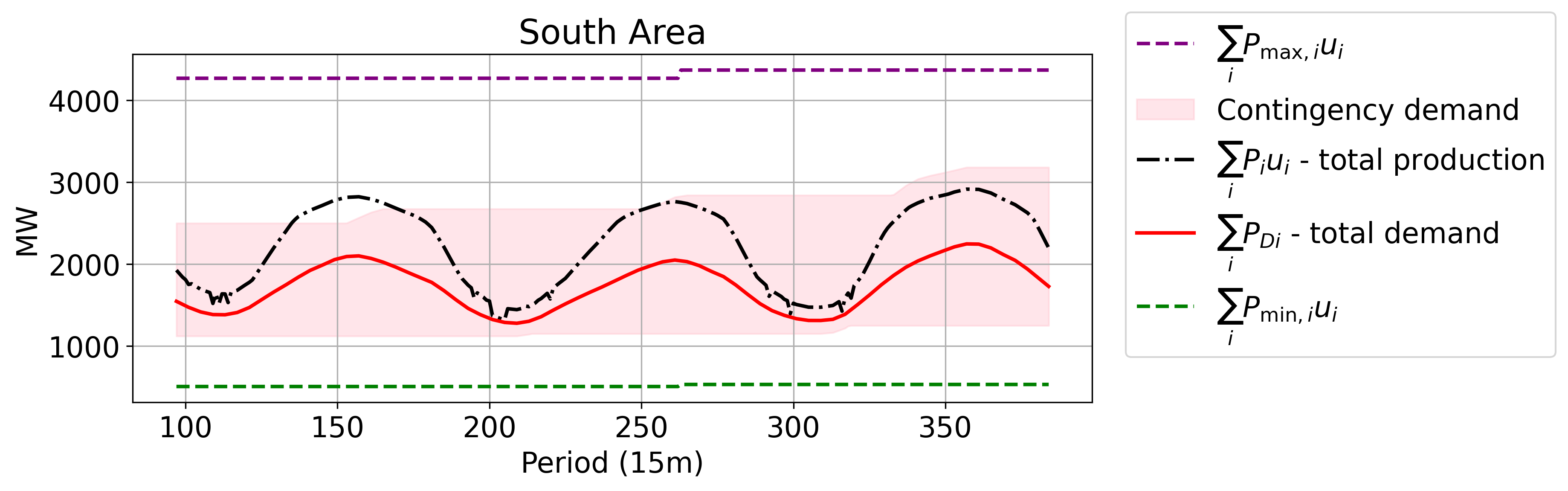}}
    \caption{Real power capacity, production, and load in the South area of the 100-bus system. The South area exports power.}
    \label{fig:south}
\end{figure*}

\begin{figure*}[!htbp]
    \centerline{\includegraphics[width=\textwidth]{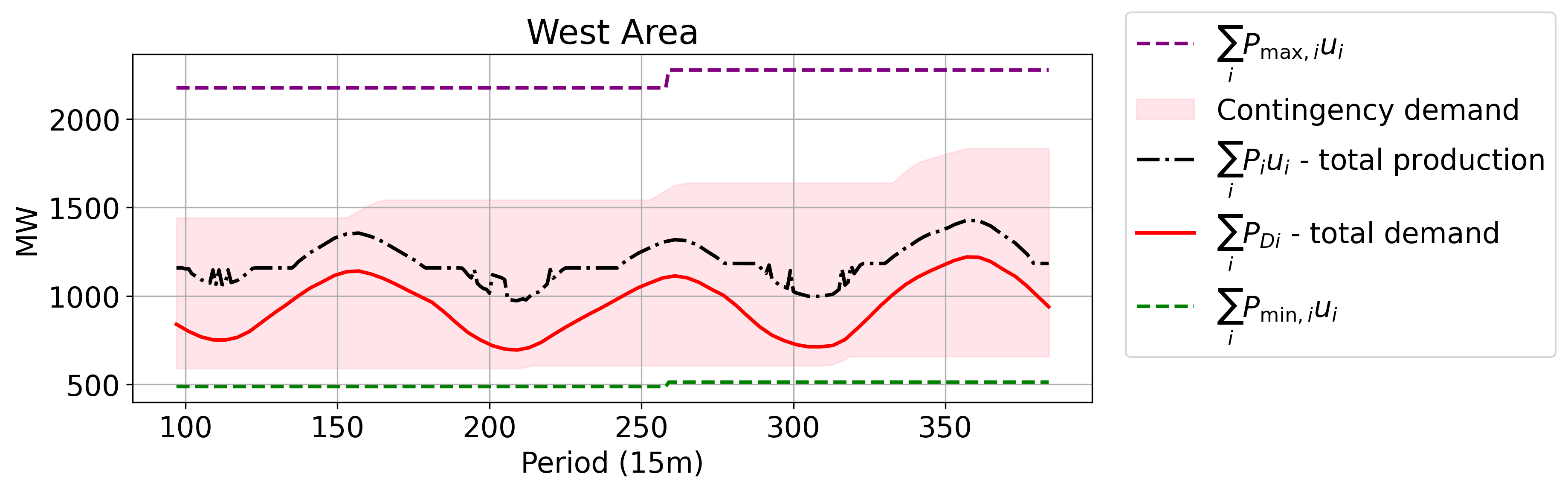}}
    \caption{Real power capacity, production, and load in the West area of the 100-bus system. The West area exports power.}
    \label{fig:west}
\end{figure*}

From the engineering viewpoint, the contingencies described in \cref{eq:ConsPDelta,eq:ConsQDelta,eq:ConsPflowdel,eq:ConsQflowdel,eq:ConsThermaldel,eq:ConsVoltagedel,eq:ConsAngledel} are meant to account for four up and down regulation contingencies and loss of a line from an interarea transmission corridor, as discussed above. \cref{fig:all areas,fig:north,fig:south,fig:west} show detailed results for the 100 bus system\footnote{We use results from the base 100-bus system to discuss these effects because the linearly scaling number of areas in the larger systems would create an unreasonable number of plots while giving only minimally more information.}. 

% We describe the resulting committed volumes in \cref{fig:all areas,fig:north,fig:south,fig:west}. For each area, the total required volume forecast is provided using the red solid line. Since the load is subject to contingency increases, the first step was to compute the corridor for load contingencies over time. The unit commitment outputs are provided as total committed \(P_{\max}\) and \(P_{\min}\) across the relevant active units as purple and green dashed lines. The total expected generation volume is shown using the dashed black line. 

\cref{fig:all areas} shows the total maximum and minimum production capabilities for on generators, total production, total load, and the load required by contingencies fluctuate over the entire 100-bus system. It shows how produced electricity (dashed-dotted black line) tracks total load (solid red line). The generation volume is subject to change after unit commitment during later market clearing procedures, but the predicted volumes offer some check on the total selected \(P_{\min}\) and \(P_{\max}\). The expected generation volume for committed units is within the bounds and is generally aligned with load forecast. This shows that \cref{alg:rruc} can meet contingencies.

\cref{fig:all areas} also shows how the committed generators change, via their changing maximum (dashed purple) and minimum (dashed green) production capabilities, when new contingencies reach the time horizon. The slight delay is due to ramp up times. New units starts to ramp up even though the currently operating set of units can meet the contingencies. This shows that \cref{alg:rruc} can change commitments to reduce costs, not just meet contingencies.

\cref{fig:north,fig:south,fig:west} shows the breakdown of \cref{fig:all areas} by area. They show that \cref{alg:rruc} is capable of finding solutions that import or export energy, even with a contingencies that hampers transmission substantially.

% The per area difference \cref{fig:north,fig:south,fig:west} offers some additional insights about how contingencies and constraints work in the presence of inter-area reserve allocation. 
 
% We use results from the base 100-bus system to discuss the area effects \cref{fig:all areas,fig:north,fig:south,fig:west}, since the linearly scaling number of zones in the larger systems would create an unreasonable number of plots while giving only minimally more information.

\cref{fig:north} shows the total maximum and minimum production capabilities for on generators, total production, total load, and the load required by contingencies within the North area. This area is connected to other areas through six lines, including two high-capacity (approximately 1200 MVA) lines 30-38 (West to North) and 68-65 (South to North). The North area imports significant amounts of energy from other areas. This makes North similar to some transit transmission operators such as Mid-continent Independent System Operator (\textbf{MISO}). The local reserve requirements are milder, partially because the North area's load would be in part supplied through the remaining interarea lines after the contingency removes a largest line. The part of \cref{fig:north} near period 340, where contingent load exceeds local area capacity, illustrates this. 

\cref{fig:south,fig:west} shows the total maximum and minimum production capabilities for on generators, total production, total load, and the load required by contingencies within the South and West areas. Both areas export substantial amounts of energy to North, explaining why their production is significantly above local load. The South area exports much more energy during peak load times because during these times the North area needs far more energy than it can supply efficiently.

\begin{figure*}[!htbp]
    \centerline{\includegraphics[width=\textwidth]{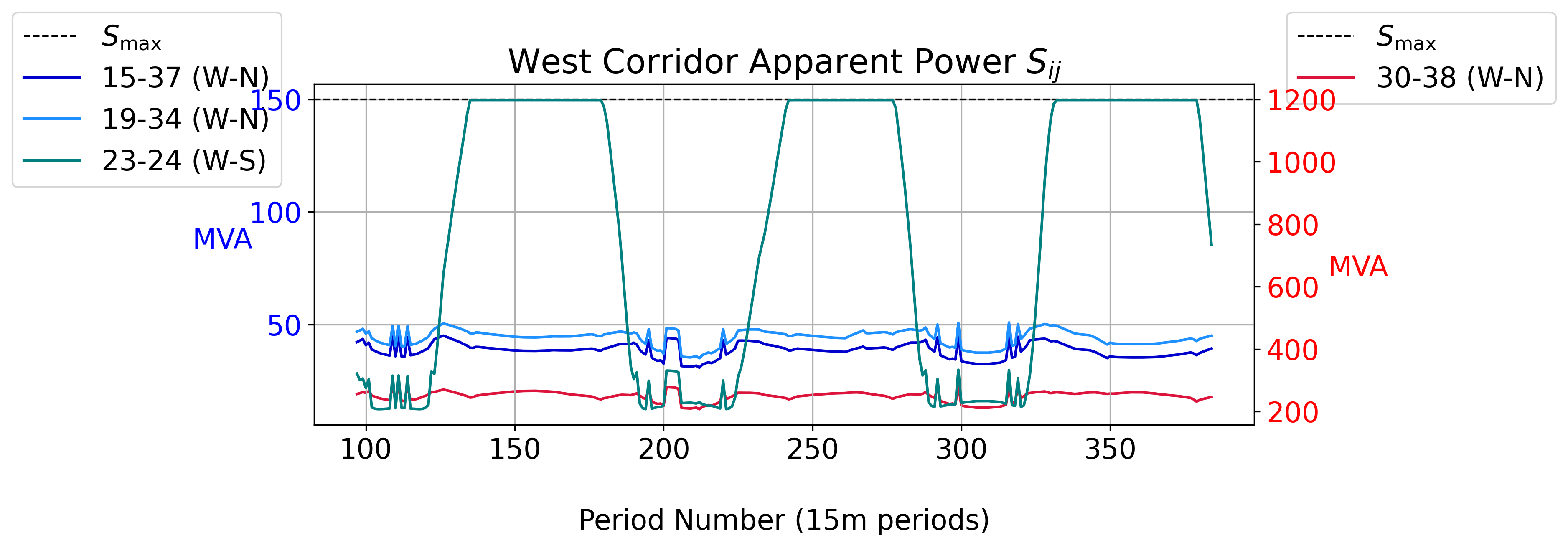}}
    \caption{Apparent power throughput (\(S = \max_{ft\in\{ij,ji\}}\sqrt{P_{ft}^2+Q_{ft}^2}\)) for lines between West and other areas. Flow between buses 23 and 24 follows daily cycles (96 periods), with minor deviations, and saturates the line during peak load. Cool colors on the left correspond to lines with an MVA limit of 150. Red on the right corresponds to a line with an MVA limit of 1200. Our contingency required UC to be feasible, even if line 30-38 goes out. This explains why it far from saturation. Lines 19-34, 23-24 retain some slack to be able transmit in case line 30-38 goes out. }
    \label{fig:west corridor s}
\end{figure*}

\begin{figure*}[!htbp]
    \centerline{\includegraphics[width=\textwidth]{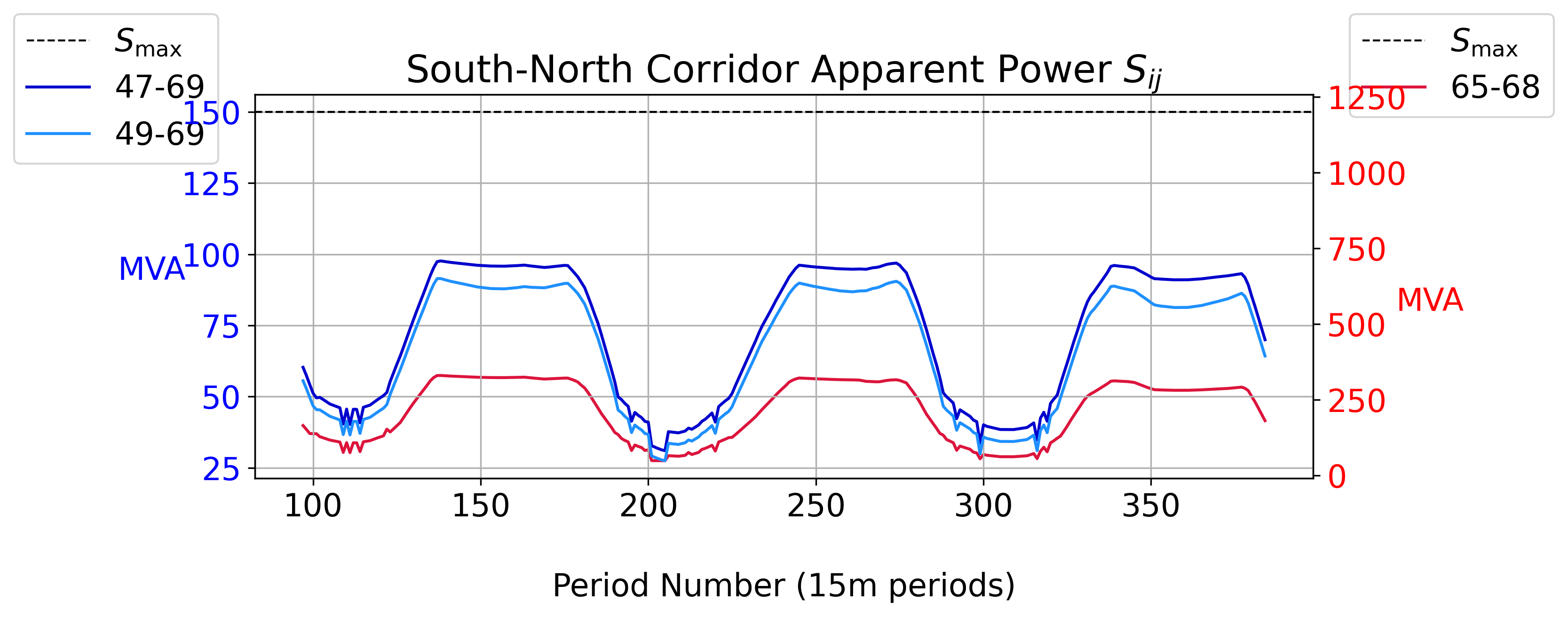}}
    \caption{Apparent power throughput (\(S = \max_{ft\in\{ij,ji\}}\sqrt{P_{ft}^2+Q_{ft}^2}\)) for lines between South and North areas. Flow between areas follows daily cycles (96 periods), with minor deviations, and is high at times of high load. Cool colors on the left correspond to lines with an MVA limit of 150. Red on the right corresponds to a line with an MVA limit of 1200. Our contingency required UC to be feasible, even if line 65-68 goes out. This explains why it far from saturation. Lines 47-69 and 49-69, retain some slack to be able transmit in case line 65-68 goes out.}
    \label{fig:north corridor s}
\end{figure*}

\cref{fig:west corridor s,fig:north corridor s} show the apparent power flow between areas. Note that the high MVA lines 30-38, and 65-68 are the furthest from saturation. This is because our contingencies required \cref{alg:rruc} to produce a solution that is feasible even if these lines go out. Note that the other lines in the corresponding corridors have enough room to absorb most or all of the apparent power on lines 30-38, and 65-68 should they go out. Any remaining difference would come from differences in generator production within the ramping limits.

\cref{fig:west corridor s} shows  power flow between the West area and other areas. Line 23-24 saturates during peak loads on daily cycles (96 periods). This shows that solutions that do not consider transmission constraints, such as in our previous work \cite{regev2026economicUC} may be infeasible. One can verify that a feasible solution works after the fact, but not recover from an infeasible one. The same is true even if one uses the UC from the problem that is not transmission constrained and then runs optimal power flow on the fully constrained system.

\cref{fig:west corridor s} shows real power flows out of the West area on the 15-37 and 19-34 lines. Real power flows into the West area on the 23-24 and 30-38 lines. 
\cref{fig:north corridor s} shows the power flow between the North and South areas. Flow between areas follows daily cycles, with small deviations. Real power flows from South to North on all lines.

%A detailed review of this pattern is shown in \cref{fig:west corridor p,fig:north corridor p}, and explains the overproduction in \cref{fig:west,fig:south}.

% By contrast, South and West are each connected to other parts of the system by four lines.
% ALL

% NORTH

% SOUTH

% WEST

% \begin{figure*}[!htbp]
%     \centerline{\includegraphics[width=\textwidth]{results_figs/west_corridor_p.png}}
%     \caption{Real power throughput (\(P_{ij}\)) for lines from West to other areas}
%     \label{fig:west corridor p}
% \end{figure*}
% \begin{figure*}[!htbp]
%     \centerline{\includegraphics[width=\textwidth]{results_figs/sn_corridor_p.png}}
%     \caption{Real power throughput (\(P_{ij}\)) for lines from South to North}
%     \label{fig:north corridor p}
% \end{figure*}


\section{Conclusions and Further Work}
\label{sec:summary}
We extended the RRUC formulation from our previous work~\cite{regev2026oneperiodUC,regev2026economicUC} to \cref{eq:UnitCommitment}, which incorporates topological constraints in a full ACOPF model. We showed experimentally that our algorithm's time complexity follows \(\mathcal{O}\left(n^{1.37}\right)\) scaling. We give experimental evidence of our algorithm maintaining accuracy and show that it can reduce production costs by optimizing over a large number of generators. Our results show that our solver can optimize a 20,000 bus system in about 7 minutes. In contrast, Gurobi cannot obtain a feasible solution to the 100 bus system in 15 minutes.

An area for future work is to explain the accuracy of the algorithm theoretically by proving an approximation bound against the optimal solution. Furthermore, one can extend the UC problem to include battery commitment and optimally shed load if it is not feasible to serve it all. 

\section*{Acknowledgements}
We thank Nicholson Koukpaizan for contributing to an improved manuscript with his thorough technical review. We thank James Nutaro for helping write the proposal that funded this work. 
\bibliographystyle{IEEEtran}
\bibliography{bibfile}

\end{document}